\documentclass[12pt]{amsart}
\usepackage{amsmath,amsfonts,amssymb,amscd,amsthm,amsbsy,epsf}
\newtheorem{prop}{Proposition}
\newtheorem{thm}[prop]{Theorem}

\def \vH{\hbox{\v H}}
\def\be{\begin{equation}}
\def\ee{\end{equation}}
\def\Cech{{\v Cech}{} }

\def \calU {{\mathcal {U}}}

\def\qed{\hbox{\hskip 6pt\vrule width6pt height7pt depth1pt \hskip1pt}}

\def\Z{{\mathbb Z}}

\def\R{{\mathbb R}}

\def\nd{\noindent}

\begin{document}

\title{Pattern-equivariant cohomology with integer coefficients}
\author{Lorenzo Sadun}
\begin{abstract}
We relate Kellendonk and Putnam's pattern-equivariant (PE) cohomology to
the inverse-limit structure of a tiling space.  This gives an
version of PE cohomology with integer coefficients,
or with values in any Abelian group.  It also provides an easy proof
of Kellendonk and Putnam's original theorem relating
PE cohomology to the \Cech cohomology of the tiling
space.  The inverse-limit structure also allows for the construction
of a new non-Abelian invariant, the PE representation variety. 
\end{abstract}
\address{Department of Mathematics, The University of Texas at Austin,
Austin, TX 78712}
\email{sadun@math.utexas.edu}
\maketitle


\markboth{Lorenzo Sadun}{P-equivariant cohomology with integer coefficients}

\medskip

\section{Background}

Kellendonk and Putnam's theory of pattern-equivariant cohomology
\cite{pattern1}
(henceforth abbreviated PE cohomology) gives a wonderfully intuitive
tool for understanding the \Cech cohomology of a tiling space.  Via
Ruelle-Sullivan currents \cite{pattern2}, it is also a powerful tool 
for understanding
the natural action of the translation group on such a space. The only
limitation is that PE cohomology is defined using
differential forms, and so necessarily has real coefficients (or
coefficients in a real vector space).  All of the additional
information contained in the integer-valued \Cech cohomology, such as
torsion and divisibility, is lost.

In this paper, we present a version of PE cohomology with
integer coefficients, or indeed with values in any Abelian group.  
This cohomology is isomorphic to the \Cech cohomology
of a tiling space with integer (or group) coefficients.  
The proof of this fact is 
quite simple, and is easily modified to give a simpler proof of Kellendonk
and Putnam's original theorem, namely that their cohomology is isomorphic
to the \Cech cohomology of the tiling space with {\em real} coefficients. 
Along the way, we also prove a de Rham theorem for branched manifolds: 
if the branched manifold is obtained by gluing polygons along edges, then
the de Rham cohomology of the branched manifold is isomorphic to the \Cech
cohomology with real coefficients.  

Finally, using the concept of
pattern-equivariance, and following an idea of Klaus Schmidt, we define a new
topological invariant of tiling spaces, called the PE representation 
variety of the tiling space with values in a group $G$.  
When $G$ is Abelian, this reduces to the first cohomology of the tiling
space with values in $G$, but when $G$ is non-Abelian, the invariant is new. 
It is related to, but different from, invariants defined by Geller and
Propp \cite{Geller-Propp} and by Schmidt \cite{Schmidt}.

Kellendonk and Putnam originally defined PE cohomology for point 
patterns.  However, tilings and point patterns are closely related. 
From any tiling, we can get a point pattern by considering the vertices; 
from any point pattern, we can get a tiling by considering Voronoi
cells.   This paper is phrased entirely in terms of tilings, and the 
edges and faces of the tilings will be essential ingredients in the 
constructions.

We consider tilings of $\R^d$ that meet three conditions: (1) there
are only a finite number of tile types, up to translation, (2) the
tiles are polygons (in $d=2$ dimensions, or polyhedra in higher
dimensions), and (3) the tiles meet full-face to full-face.  The
second and third conditions may at first appear severe, but any tiling
that has finite local complexity (FLC) with respect to translation is
mutually locally derivable (MLD) to a tiling that meets these
conditions.  A tiling with FLC is associated with a Delone point pattern
with FLC, which in turn is associated with a derived Voronoi tiling that 
meets the three conditions.  For details of this construction, see
\cite{Frank}.

The translation group $\R^d$ acts naturally on a tiling.  If $x\in \R^d$
and $T$ is a tiling, then $T-x$ is the tiling $T$ translated backwards by $x$.
That is, a neigborhood of the origin in $T-x$ looks like a neighborhood of
$x$ in $T$.  Two tilings are considered $\epsilon$-close if they
agree, up to a rigid translation by $\epsilon$ or less, on a ball of 
radius $1/\epsilon$ around the origin.  Given a tiling $T$, the {\em 
tiling space} of $T$, or the {\em continuous hull} of $T$, denoted 
$\Omega_T$, is the closure (in the space of all tilings) of the translational
orbit of $T$.  Another characterization of $\Omega_T$ is the following: 
a tiling $T'$ is in $\Omega_T$ if and only if every patch of $T'$ is 
found (up to translation) somewhere in $T$.   Since $T$ has finite local
complexity, $\Omega_T$ is compact.  If $T$ is {\em repetitive}, meaning 
that every pattern in $T$ repeats with bounded gaps, then $\Omega_T$ is 
a minimal dynamical system, i.e., every translational orbit in $\Omega_T$ 
is dense.  

Let $T$ be a tiling of $\R^d$.  Two points $x,y \in \R^d$ are 
{\em $T$-equivalent
to radius $r$} if $T-x$ and $T-y$ agree exactly on a ball of radius $r$ around
the origin.  A function $f$ on $\R^d$ is 
{\em PE with radius $r$} with respect to $T$ if there
exists an $r>0$ such that $f(x)=f(y)$
whenever $x$ and $y$ are $T$-equivalent to radius $r$. A differential $k$-form
on $\R^d$ is PE if it can be written as a finite sum
$\sum_I f_I(x) dx^I$, where each $f_I(x)$ is a PE
function and $dx^I = dx^{i_1}\wedge \cdots \wedge dx^{i_k}$ is a 
constant $k$-form.  (In \cite{pattern2} a distinction is made between 
functions that meet the above conditions, which are called {\em strongly}
PE, and {\em weakly} PE
functions that are uniform limits of strongly PE functions.  In this paper
we will only consider strongly PE functions.)

Let $\Lambda^k_P(T)$ denote the set of PE $k$-forms.  It is easy to see 
that the exterior derivative of a PE form is PE, and $d^2=0$, as usual,
so we have a differential complex:
\be 0 \longrightarrow\Lambda^0_P(T)
\xrightarrow{d_0} \Lambda^1_P(T) \xrightarrow{d_1}
\cdots \xrightarrow{d_{n-1}}
\Lambda^d_P(T) \xrightarrow{d_d} 0. 
\ee
The pattern-equivariant cohomology of $T$ is the cohomology of this complex:
\be H^k_P(T)= (\ker d_k)/(\text{Im}(d_{k-1})). \ee
Kellendonk and Putnam proved:
\begin{thm} \label{real} $H^*_P(T)$ is isomorphic to $\vH^*(\Omega_T, \R)$,
the \Cech cohomology with real coefficients of the tiling space of $T$.
\end{thm}

In section 2 we review the structure of the tiling space 
$\Omega_T$ as an inverse limit of approximants $\Gamma_k$, and give a 
new proof of Kellendonk and Putnam's theorem.  
We will see that a PE form of $\R^d$ is just the pullback of a form on
$\Gamma_k$, and that the PE cohomology of $T$ is just the direct limit 
of the de Rham cohomology of the approximants. Thanks to a de Rham theorem
for branched manifolds (proved in the appendix), this is the same as the
direct limit of the \Cech cohomology of the approximants with real 
coefficients, and hence equal to the \Cech cohomology of the inverse limit. 

In section 3 we use the CW structure of $\R^d$ induced by the tiling $T$ to
define PE cellular cochains with values in an arbitrary Abelian group.  
Repeating the argument of section 2, these are the pullbacks of cellular
cochains on approximants $\Gamma_k$, so the PE cohomology of $T$ is 
isomorphic to the
direct limit of the cellular cohomology of $\Gamma_k$, hence to the 
direct limit of the \Cech cohomology of $\Gamma_k$, hence to the \Cech
cohomology of $\Omega_T$.  

In section 4 we define a new topological invariant, the PE
representation variety of a tiling with gauge group $G$.  
This is defined intrinsically in terms
of PE connections and PE gauge transformations, but is seen to equal the
direct limit $\lim_\to Hom(\pi_1(\Gamma_n), G)/G$.  For each approximant, 
the representation variety is derived from the fundamental group, but in 
the limit the fundamental group disappears, while the representation variety
remains.

\section{Inverse limits and real cohomology}

There are several constructions of tiling spaces as inverse limits 
\cite{AP, ors, bbg, gambaudo, Sadun}.  
For our purposes, the most useful is G\"ahler's construction, 
building on the earlier work of Anderson and Putnam \cite{AP}.  
We sketch the construction below; see \cite{Gahler, Sadun} for 
further details and generalization. 

A point in the $n$-th
approximant $\Gamma_n$ is a description of a tile containing the origin,
its nearest neighbors (sometimes called the ``first corona''), its
second nearest neighbors (the ``second corona'') and so on out to the
$n$-th nearest neighbors.  
The map $\Gamma_n \to \Gamma_{n-1}$ simply
forgets the $n$-th corona.  A point in the inverse limit is then a
consistent prescription for constructing a tiling out to infinity.  In
other words, it defines a tiling.

What remains is to construct $\Gamma_n$ out of geometric pieces. 
We consider two tiles $t_1$, $t_2$ in $T$ to be
equivalent if a patch of $T$, containing $t_1$ and its 
first $n$ coronas, is identical, up to translation, to a 
similar patch around $t_2$.   Since $T$
has finite local complexity, there are only finitely many
equivalence classes, each of which is called an $n$-collared tile. 
(Note that an $n$-collared tile is no bigger than an ordinary tile,
but its label contains information on how the first $n$ coronas are 
laid out around it.)

For each $n$-collared tile $t_i$, we consider how such a tile
can be placed around the origin.  This is tantamount to picking a point
in $t_i$ to call the origin, so the set of ways to place $t_i$ is
just a copy of $t_i$ itself. 

A patch of a tiling in which the origin is on the boundary of two or more tiles
is described by points on the boundary of two or more $n$-collared tiles,
and these points must be identified.  The branched manifold $\Gamma_n$ is 
the disjoint union of the $n$-collared tiles $t_i$, 
modulo this identification. 
Each of the points being identified
carries complete information about the placement of all the tiles that meet
the origin, together with their first $n-1$ coronas. 

$\Gamma_n$ is called the {\em Anderson-Putnam complex} with $n$-collared tiles.
This is a branched manifold, with branches where multiple tile boundaries 
are identified.  At the branches, there is a well-defined tangent space, 
with a global trivialization given by the action of the translation group, so
it makes sense to speak of smooth functions and smooth differential forms.

The tiling $T$ induces a natural map $\pi_n: \R^d \to \Gamma_n$. 
For each $x \in \R^d$, $\pi_n(x)$ is the point in $\Gamma_n$ that describes
a neighborhood of the origin in $T-x$, or equivalently a neighborhood of $x$
in $T$.  Let $L$ be the diameter of the largest tile.  If $x$ and $y$ are
$T$-equivalent to radius $R$, and if $R>(n+1)L$, then $\pi_n(x)=\pi_n(y)$. 
Likewise, if $\pi_n(x)=\pi_n(y)$, then there is a radius $r$ (which grows
uniformly with $n$) such that $x$ and $y$ are $T$-equivalent to radius $r$. 
Thus a function $f$ of $\R^d$ is PE if and only if there exists an $n$
such that $f(x)=f(y)$ whenever $\pi_n(x)=\pi_n(y)$.  That is, there exists
a function on $\Gamma_n$ whose pullback is $f$, in which case there also
exist functions on all $\Gamma_{n'}$ with $n'>n$ whose pullback is $f$. 
The same argument applies to differential forms.  That is:
\begin{thm} \label{pullback} 
$\Lambda^k_P(T) = \cup_n \pi_n^*(\Lambda^k(\Gamma_n))$,
\end{thm}
\noindent
where $\Lambda^k(\Gamma_n)$ denotes the set of smooth $k$-forms on $\Gamma_n$.

Our last ingredient for proving Theorem (\ref{real}) is the de Rham theorem
for branched manifolds:
\begin{thm} \label{deRham}  
  If $X$ is a branched manifold obtained by gluing tiles along their
  common boundaries, then the de Rham cohomology of $X$ is naturally
  isomorphic to the \Cech cohomology of $X$ with real coefficients.
\end{thm}

\nd Proof:  See appendix. 

\nd Proof of Theorem (\ref{real}): Every PE cohomology class is
represented by a closed PE form $\omega$, which is the pullback of a
closed form $\omega^{(n)}$ on $\Gamma_n$, which defines a de Rham
cohomology class $[\omega^{(n)}]_{dR}$ on $\Gamma_n$, which defines a
\Cech cohomology class on $\Gamma_n$.  This defines a class in the
direct limit of the \Cech cohomologies of the approximants, which is
the \Cech cohomology of the inverse limit space $\Omega_T$.  In other
words, we have a map $\phi: H_P^*(T) \to \vH^*(\Omega_T)$.  We must
show that $\phi$ is well-defined, 1--1 and onto.

Being well-defined means being independent of the choices made, namely
of $n$ and of the representative $\omega$.  If we pick a different
approximant $\Gamma_{n'}$, with $n'>n$, then there are commuting maps:
\begin{center}
\begin{picture}(100,80)
\put(47,70){$\Omega_T$}   
\put(14,46){$\pi_{n'}$}
\put(44,60){\vector(-3,-4){26}}  
\put(64,60){\vector(2,-3){24}}   
\put(82,46){$\pi_n$}
\put(0,10){$\Gamma_{n'}$}
\put(90,10){$\Gamma_n$}
\put(22,14){\vector(3,0){60}}    
\put(46,6){$\rho_{nn'}$}
\end{picture}
\end{center} 
The form $\omega^{(n')}$ is the pullback by the forgetful map
$\rho_{nn'}$ of the form $\omega^{(n)}$. By
the naturality of the isomorphism between de Rham and \Cech cohomology, 
the \Cech class defined by $\omega^{(n')}$ is the pullback of the \Cech
class defined by $\omega^{(n)}$, and hence defines the same element of 
the direct limit of \Cech cohomologies. 

Next suppose that $\omega'$ and $\omega$ define the same PE class of
$T$.  Then $\omega' = \omega + d\nu$, where $\nu$ is a PE form.  Then
there exists an $n$ such that $\omega$, $\nu$ and $\omega'$ are all
pullbacks of forms on $\Gamma_n$.  But then $\omega'{}^{(n)} =
\omega^{(n)} + d \nu^{(n)}$, so $\omega$ and $\omega'$ define the same
de Rham cohomology class on $\Gamma_n$, hence the same \Cech class,
and hence the same class in the direct limit of cohomologies.  This shows that
$\phi$ is well-defined. 

To see that $\phi$ is 1--1, suppose that $\phi([\omega])=0$.  Then there
exists a finite $n$ such that the \Cech class defined by $\omega^{(n)}$ is 
zero, so the de Rham class defined by $\omega^{(n)}$ is zero, so there 
is a form $\nu$ on $\Gamma_n$ with $d\nu = \omega^{(n)}$, so 
$[\omega] = [d \pi_n^*(\nu)]=0$ in PE cohomology. 

Finally, every class $\gamma$ 
in the direct limit of \Cech cohomologies is represented
by a \Cech class in some $\Gamma_n$, which in turn is represented by a closed
form on $\Gamma_n$.  Pulling this back to $\R^d$, we get a closed 
PE form, hence a PE class whose image under $\phi$ is $\gamma$.  \qed

\section{Integer-valued PE cohomology}
  
We next construct PE cohomology with integer coefficients. 
Although we speak of integers, the same arguments apply, with no loss of 
generality, to coefficients in any Abelian group. 

A tiling is more than a point pattern.  A tiling $T$ gives a CW decomposition
of $\R^d$ into 0-cells (vertices), 1-cells (edges), 2-cells (faces), and so
on.  We may therefore speak of cellular chains with the obvious boundary
maps, and we can dualize to get cellular cochains with integer (or group)
coefficients.   Likewise, each approximant $\Gamma_n$ is constructed as a
union of tiles, glued along boundaries.  $\Gamma_n$ is a CW complex, and the
map $\pi_n: \R^d \to \Gamma_n$ is cellular. 

For $k>0$, we say that two $k$-cells $c_1$ and $c_2$ are {\em
$T$-equivalent to radius $r$} if there exist interior points $x \in c_1$
and $y\in c_2$ such that $x$ and $y$ are $T$-equivalent to radius $r$.
($T$-equivalence for 0-cells has already been defined.) We say a cochain
$\alpha$ on $\R^d$ is PE if there exists an $r>0$ such that
$\alpha(c_1)= \alpha(c_2)$ whenever $c_1$ and $c_2$ are equivalent to
radius $r$.  The coboundary $\delta$ of a PE cochain is PE (albeit with a 
slightly larger radius).  Denoting the PE $k$-cochains on $\R^d$ by
$C^k_P(T)$, we have a complex
\be 0 \longrightarrow C^0_P(T)
\xrightarrow{\delta_0} C^1_P(T)
\xrightarrow{\delta_1} \cdots \xrightarrow{\delta_{n-1}} 
C^d_P(T) \xrightarrow{\delta_d} 0. 
\ee
and define the integer-valued PE cohomology of $T$ to be the cohomology
of this complex. 

\begin{thm} \label{integer}
The integer-valued PE cohomology of $T$ is isomorphic to the \Cech cohomology
of $\Omega_T$. 
\end{thm}

\nd Proof:   As before, we see that a cochain on $\R^d$ is PE if and
only if it is the pullback of a cochain on $\Gamma_n$ for some $n$ (and 
hence any sufficiently large $n$). 
Every PE cohomology class is 
represented by a closed PE cochain $\gamma$, which is the pullback of 
a closed cochain $\gamma^{(n)}$ on $\Gamma_n$, which defines a cellular
cohomology class $[\gamma^{(n)}]$ on $\Gamma_n$, which defines a 
\Cech cohomology class on $\Gamma_n$, which defines a class in the 
direct limit of the \Cech cohomologies of the approximants, which is 
the \Cech cohomology of the inverse limit space $\Omega_T$.  

Showing that this compound map from $H^*_P(T)$ to $\vH^*(\Omega_T)$
is well-defined and an isomorphism proceeds exactly as in the
proof of Theorem (\ref{real}), with the
word ``form'' replaced by ``cochain'', ``exterior derivative'' 
replaced by ``coboundary'', and ``de Rham theorem'' replaced by ``the natural
isomorphism between cellular and \Cech cohomology'', which holds for all 
CW complexes.  \qed

\section{PE representation varieties}

Gauge theory provides a wealth of topological invariants of manifolds.
Given a manifold $M$ and a connected Lie Group $G$, we can consider
the moduli space of flat connections on the (trivial) principal bundle
$M\times G$, modulo gauge transformation.  At first glance this moduli
space might appear to be only an invariant of the smooth structure of
$M$, but in fact it is a topological invariant.  Flat connections are
described by their holonomy, i.e., by a homomorphism $\pi_1(M) \to G$,
and gauge equivalent connections have maps that are related by
conjugation by an element of $G$.  In other words, the moduli space is
isomorphic to the {\em representation variety} $Hom(\pi_1(M), G)/G$,
where the quotient is by conjugation.  (The representation variety is
well-defined even when $G$ is not connected, and is a non-Abelian
generalization of the first cohomology of $M$.  When $G$ is Abelian,
$Hom(\pi_1(M), G)/G = Hom(\pi_1(M), G)= Hom(H_1(M),G) = H^1(M,G).$)

Using pattern equivariance, we can play the same game for tilings.
Given a tiling $T$ and a connected Lie group $G$, we consider flat PE
connections on the trivial principal bundle $\R^d \times G$, modulo PE
gauge transformations.  We call this the {\em PE representation
  variety} of the tiling $T$.  As with PE cohomology, the PE
representation variety is defined abstractly, without reference to
inverse limits, but is most easily understood in terms of inverse
limits and G\"ahler's construction.

If $A$ is a PE connection 1-form with zero curvature, then $A$ is the 
pullback of a connection 1-form on $\Gamma_n$ with zero curvature, and so
defines an element of the representation variety 
$Hom(\pi_1(\Gamma_n), G)/G$.
Fundamental groups do not behave well under inverse limits, but their
duals do.  For $n' > n$, we have a forgetful map $\Gamma_{n'} \to \Gamma_n$,
which induces a map $\pi_1(\Gamma_{n'}) \to \pi_1(\Gamma_n)$, which induces
a map $Hom(\pi_1(\Gamma_{n}),G) \to Hom(\pi_1(\Gamma_{n'}), G)$.  This commutes
with the $G$ action, and induces a map of representation varieties:
$Hom(\pi_1(\Gamma_{n}),G)/G \to Hom(\pi_1(\Gamma_{n'}), G)/G$.  
The PE representation variety of $T$ is the direct limit
$\lim_{\rightarrow} Hom(\pi_1(\Gamma_{n}),G)/G$. 

Notice that for each approximant $\Gamma_n$, the representation
variety contained no additional information to that in the fundamental
group of $\Gamma_n$.  However, this changes in the infinite limit.
The fundamental group of $\Omega_T$ is trivial, but the PE
representation variety of $T$ can be nontrivial.  Indeed, when $G$ is
Abelian, the PE representation variety of $T$ is the first \Cech
cohomology of $\Omega_T$ with coefficients in $G$.

It is natural to ask whether there is a limiting object, dual to the PE 
representation
variety, that plays a role analagous to the fundamental group.  
For $\Z^d$ actions on subshifts of finite type, something akin to this was 
considered by Geller and Propp \cite{Geller-Propp}.  
Their tiling group is somewhat smaller than 
$\pi_1(\Gamma_n)$, as they only consider paths that correspond to closed loops
in $\R^d$.  That is, they consider the kernel of the map $\pi_1(\Gamma_1)
\to \R^d$, $\gamma \to \int_\gamma d \vec x$.  They then consider the
{\em projective fundamental group} of a tiling to be the direct limit of 
the tiling groups of the approximants.  The cohomology class of Schmidt's {\em
fundamental cocycle} \cite{Schmidt} is a non-Abelian structure dual to
Geller and Propp's tiling group.  As such, it should be closely
related to the PE representation variety, but this relation is not yet
understood.

\section{Acknowledgements}

Many thanks to Jean Bellissard, Johannes Kellendonk, Ian Putnam, Bob
Williams and especially Klaus Schmidt for helpful discussions.  In
particular, the idea for the PE representation variety grew directly
out of discussions with Schmidt.  This work is partially supported by
the National Science Foundation. 

\section{Appendix: A de Rham theorem for branched manifolds}

Let $X$ be a branched manifold obtained by gluing Euclidean
polygons (or polyhedra) along their common faces, as in the
Anderson-Putnam complex.  Note that the branch set is locally
star-shaped with respect to every point on the branch set.  That is,
for each point $p$, there is a ball $U_p$ around $p$ such that
the branch set, intersected with $U_p$, is a union of (Euclidean)
line segments with endpoint $p$. If $p$ is a vertex, we 
take $U_p$ to be a ball of radius less half than the distance from $p$
to the nearest other vertex.  (This ensures that the neighborhoods
corresponding to different vertices do not intersect, a condition that
will prove useful.) If $p$ is on an edge, we take the radius
of $U_p$ to be less than half the distance to the nearest other edge,
and also less than the distance to the nearest vertex.  If $p$ is on a
face, we take the radius to be less than half the distance to the
nearest other face, and also less than the distance to the nearest
edge or vertex.

\begin{thm}[Poincar\'e Lemma] \label{Poincare}
Let $U = U_{p_1} \cap U_{p_2} \cap \cdots \cap U_{p_n}$ be a non-empty
intersection of a finite number of balls $U_{p_i}$ constructed as
above. Then $U$ has the de Rham cohomology of a point.  That is,
$H^0_{dR}(U)=\R$, and $H_{dR}^k(U)=0$ for $k>0$. 
\end{thm}

\nd {\bf Proof:} If $U$ does not touch the branch set of $X$, then $U$
is a convex subset of a single polygon, and the usual proof of the
Poincar\'e Lemma applies.  We therefore assume that $U$ hits the
branch set of $X$.  But by the construction of $U_{p_i}$, this means
that all the points $p_i$ must lie on
the same $k$-cell $C$, and that the branch set, intersected with $U$,
is $U \cap C$. Since the branch set of each $U_{p_i}$ is a convex ball in $C$,
the branch set of $U$ is a convex set in $C$, and in particular is
star-shaped. Let $p_0$ be any point in $U\cap C$. 

The proof of the Poincar\'e Lemma found in \cite{Spivak}, 
page 94, then carries over almost word for word. We define a homotopy operator
$I$ that takes $k$-forms to $k-1$ forms such that, for any $k$-form $\omega$,
$d(I\omega) + I (d\omega)=\omega$.  If $\omega$ is closed, then $\omega = 
d I \omega$ is exact.  The form $I \omega$ is computed by integrating 
$\omega$ out along a straight line from $p_0$.  Since the straight line
is either entirely in the branch set or entirely outside the branch set, 
the value of $I(\omega)$ at a point does not depend on which coordinate disk
is used to do the calculation.  Specifically, 
\begin{equation}
(I \omega)_{i_1,\ldots,i_{k-1}}(x) = \sum_j \int_0^1 t^{k-1}(x-p_0)_j 
\omega_{j,i_1, i_2,\ldots,i_{k-1}}(p_0 + t(x-p_0)) dt. \qed
\end{equation}

Now take an open cover ${\calU}$ of $X$, with each open set of the form
$U_{p_i}$. 

\begin{thm} \label{Cech-cover}
The de Rham cohomology of $X$ is isomorphic to the Cech cohomology of
$\calU$ with real coefficients. 
\end{thm}

\nd {\bf Proof}:  We set up the 
\Cech-de Rham double complex as in Chapter II of \cite{Bott-Tu}. The rows are
all exact, by the partition-of-unity argument found in (\cite{Bott-Tu}, 
Proposition 8.5).  This implies that the cohomology of the double complex
is isomorphic to the de Rham cohomology of $X$.  Likewise, by the
Poincar\'e Lemma, the columns are all exact, and the cohomology of the double
complex is isomorphic to the \Cech cohomology of the cover with 
real coefficients. \qed

\nd {\bf Proof of Theorem \ref{deRham}}:  The Cech cohomology of $X$
is, by definition, the direct limit of the Cech cohomology of the open
covers of $X$.  However, every open cover has a refinement of the form
$\calU$, since each set $U_{p_i}$ is made by taking an arbitrarily
small ball around the point $p_i$. In computing the Cech cohomology of
$X$, it is thus sufficient to consider only covers of the form
$\calU$. For each of these, theorem \ref{Cech-cover} gives an isomorphism
to de Rham.  So each term in the direct limit is the same, and each
map is an isomorphism. \qed

\nd {\bf Remark:} The proof given here can easily be generalized to
broader classes of branched manifolds.  The only necessity is for the
branch set to be sufficiently well-behaved that one can prove a
Poincare Lemma.  An arbitrary branched surface may not have a
well-behaved branch set, but Williams \cite{Williams} showed how to construct
inverse limit spaces using only ``nice'' branched manifolds, for which
a theorem similar to \ref{deRham} should be expected to hold.

\end{document}